\documentclass[12pt]{article}
\usepackage[cp1251]{inputenc}
\usepackage[russian]{babel}
\usepackage{amssymb,amsmath}
\usepackage[dvips]{graphicx,color}
\usepackage{mathrsfs}

\usepackage{amssymb,amsmath}
\usepackage{mathrsfs}
\usepackage{amssymb}
\usepackage{amsthm}

\newtheorem{Theorem}{Theorem}

\newtheorem{Definition}{Definition}

\date{}

\begin{document}

\title{On stability of triangular factorization of positive operators}
\author{Mikhail I. Belishev
\thanks{St. Petersburg Department of Steklov Mathematical Institute
of Russian Academy of Sciences, Fontanka 27, St. Petersburg,
Russia, 191023; belishev@pdmi.ras.ru}, Aleksei F.
Vakulenko\thanks{St. Petersburg Department of Steklov Mathematical
Institute of Russian Academy of Sciences, Fontanka 27, St.
Petersburg, Russia, 191023; vak@pdmi.ras.ru}}

\maketitle

\begin{abstract}
Let $\mathfrak f=\{\mathscr F_s\}_{s>0}$ be a nest and $C$ a
bounded positive operator in a Hilbert space $\mathscr F$. The
representation $C=V^*V$ provided $V\mathscr F_s\subset\mathscr
F_s$ is a triangular factorization (TF) of $C$ w.r.t. $\mathfrak
f$. The factorization is stable if
$C^\alpha\underset{\alpha\to\infty}\to C$ and
$C^\alpha=V^{\alpha\,*}V^\alpha$ implies $V^\alpha\to V$. If $C$
is positive definite (isomorphism), then TF is stable. The paper
deals with the case of positive but not positive definite $C$. We
impose some assumptions on $C^\alpha$ and $C$ which provide the
stability of TF.
\end{abstract}

\section{Introduction}

\noindent$\bullet$\,\,\, Triangular factorization (TF) of
operators is a classic tool for solving inverse problems of
analysis and mathematical physics (M.G.Krein, V.A.Mar\-chenko,
I.M.Gelfand, B.M.Levitan, L.D.Faddeev, R.Newton et al). The
Boundary Control method \cite{B Obzor IP 97,B UMN,B IPI} also uses
TF. The method (more precisely, its time-domain version) provides
time-optimal reconstruction of parameters, manifolds and metrics.
Until recently, there were no results on its stability, that is,
the continuity of the correspondence between inverse data and the
reconstructed objects. A possible approach to its study is
analysis of TF stability. This path faces significant obstacles in
the case when the factorized operator is not an isomorphism. In
applications, the latter is typical for multidimensional inverse
problems and reflects their strong ill-posedness.

Though the motivation comes from concrete inverse problems, we
study the stability of TF at general operator level, without
referring to applications. At the same time, we hope that our
results and their further development will be useful to them.
\smallskip

\noindent$\bullet$\,\,\,Unless otherwise stated, we deal with
bounded operators.

Let $\mathfrak f=\{\mathscr F_s\}_{s\in[0,T]}$ be a {\it nest}  in
a Hilbert space $\mathscr F$, i.e., a family of extending
subspaces \cite{GK,Dav} \footnote{in the Russian literature the
term {\it chain} is in the use \cite{GK}}. Let $X_s$ be
(orthogonal) projection in $\mathscr F$ onto $\mathscr F_s$. For
future applications purposes, we assume $\mathfrak f$ to be
continuous and bordered, which means
\begin{equation}\label{Eq def nest}
\mathscr F_s\subset\mathscr F_{s'},\,\,\,0\leqslant s<s'\leqslant
T;\quad \text{s\,-}\!\lim\limits_{s'\to s}X_{s'}=X_s;\quad
\mathscr F=\{0\},\,\,\,\mathscr F_T=\mathscr F.
\end{equation}

Let $C$ be a positive operator in $\mathscr F$, i.e., $(Cf,f)> 0$
holds for $f\not=0$. The representation
$$
C=V^*V\quad\text{provided}\quad VX_s=X_sVX_s,\quad s\in[0,T]
$$
(equivalently $V\mathscr F_s\subset\mathscr F_s$ ) is a {\it
triangular factorization} (TF) of $C$ with respect to the nest
$\mathfrak f$. The factorization is said to be stable if
$C^\alpha\underset{\alpha\to \infty}\to C$ and
$C^\alpha=V^{\alpha\,*}V^\alpha$ implies $V^\alpha\to V$ (with a
relevant understanding of the convergences). If $C$ is of the form
$C=\Bbb I+K$ with a compact $K$, and is positive definite, i.e.,
$(Cf,f)\geqslant\gamma\|f\|^2$ with $\gamma>0$ holds, then TF is
realized by the classical construction of the triangular
truncation operator integral (M.S.Brodskii, I.Ts.Gokhberg and
M.G.Krein \cite{GK}), which provides a stable (so-called {\it
special}) factorization $C=(\Bbb I+L)^*(\Bbb I+L)$ \cite{GK,Dav}.

Meanwhile, there is another construction ({\it amplitude
integral}, AI), ge\-ne\-ralizing the Brodskii construction to a
wider class of factorizable operators and providing their
canonical TF \cite{B POMI PREP_90,BKach Oper Integral 94}. For
factorization by AI, the operator does not have to be isomorphic
(positive definite). However, the question arises whether the
factorization by AI is stable. We show that, under certain
additional assumptions about $C^\alpha$ and $C$, the answer is
affirmative.
\smallskip

\noindent$\bullet$\,\,\, The key element of AI is the operator
{\it diagonal}. For an operator $W$ and a nest $\mathfrak f$, it
is defined by
\begin{equation}\label{Eq def Diag}
D_W={\rm w}\text{-}\lim\limits_{r\to
0}\,\sum\limits_{k=0}^N(P_{s_k}-P_{s_{k-1}})\,W\,(X_{s_k}-X_{s_{k-1}})=:\int_{[0,T]}dP_s\,W\,dX_s,
\end{equation}
where $0=s_0<s_1<\dots<s_N=T$ is a partition of $[0,T]$ of the
range $r$, $X_s$ and $P_s$ project onto $\mathscr F_s$ and
$\overline{W\mathscr F_s}$ respectively. The  limit, if it exists,
is called the diagonal of $W$ w.r.t. the nest $\mathfrak f$.
Diagonal intertwines the projections:
\begin{equation}\label{Eq DX=PD}
D_WX_s=P_sD_W,\qquad D^*_WP_s=X_sD^*_W
\end{equation}
holds.

We say that a positive operator $C$ admits a {\it canonical} TF
w.r.t. a nest $\mathfrak f$ if its (positive) square root has the
diagonal $D_{\sqrt C}$, which satisfies $D_{\sqrt C}\,D^*_{\sqrt
C}=\Bbb I$, and the representation
\begin{equation}\label{Eq 0}
C=V^*V,\qquad  V=D^*_{\sqrt C}\,{\sqrt C}
\end{equation}
holds with the factor $V$ satisfying $VX_s=X_sVX_s$ by (\ref{Eq
DX=PD}). Such a notion is motivated by applications \cite{B Obzor
IP 97,B IPI,BSim_AA_2024,BSim_JOT_2025}.

Our main result is Theorem \ref{T0}, which establishes the
stability of the canonical TF as follows.
\begin{Definition}\label{D reg conv}
In the notation from (\ref{Eq def Diag}), we say that $W^\alpha$
converges to $W$ {\it regularly on the nest} $\mathfrak f$  if the
relations ${\rm s}\text{-}\!\lim\limits_{\alpha\to
\infty}W^\alpha=W$ and ${\rm s}\text{-}\!\lim\limits_{\alpha\to
\infty}P^\alpha_s= P_s$ hold for all $s\in[0,T]$. In this case we
write $W^\alpha\overset{\rm reg}\to W$.
\end{Definition}
At the end of the article, in the "Comments" section, the examples
are given which support and motivate the introduction of this kind
of convergence.
\begin{Theorem}\label{T0}
Let $C^\alpha\overset{\rm s}\to C$, operators $C$ and $C^\alpha$
admit the canonical TF (\ref{Eq 0}), and
$\sqrt{C^\alpha}\overset{\rm reg}\to \sqrt C$ holds. Let the
integral sums, which define  $D_{\sqrt C}$ and $D_{\sqrt
{C^\alpha}}$ by (\ref{Eq def Diag}), converge to their (weak)
limits uniformly w.r.t. to $\alpha$. Then the convergence of the
triangular factors $V^\alpha\overset{\rm w}\to V$ occurs.
\end{Theorem}

\noindent Its proof is preceded by relevant preliminaries. We
emphasize that $C$ and $C^\alpha$ are assumed to be positive, but
not necessarily positive definite.
\smallskip

\noindent$\bullet$\,\,\, At the heuristic level, the AI and
diagonal were introduced in \cite{B POMI PREP_90}. Later on its
construction was rigorously justified and developed \cite{BKach
Oper Integral 94,B Obzor IP 97,BSim_AA_2024}; its relationship
with the classical triangular truncation integral has been
recognized \cite{B TF Isom,BPush} \footnote{Paper \cite{BPush},
due to the fault of the first author, contains an erroneous
statement (Proposition 2.6)}.

The authors express their gratitude to S.A. Simonov for valuable
discussions of the results of the article.

\section{Canonical TF and its stability}

\subsubsection*{Triangularity}

So, we deal with a Hilbert space $\mathscr F$ and a family of its
(closed) subspaces $\mathfrak f=\{{\mathscr F}_s\}_{s\in[0,T]}$,
which is a {\it nest}, i.e., satisfies (\ref{Eq def nest}). An
operator $V$ in $\mathscr F$ is {\it triangular} (with respect to
the nest $\mathfrak f$) if $V\mathscr F_s\subset\mathscr F_s$ is
valid for all $s$. The latter is equivalent to $VX_s=X_sVX_s$.

An operator $C$ admits the triangular factorization (w.r.t.
$\mathfrak f$) if there is a triangular operator $V$ providing
$C=V^*V$. This of course implies $C\geqslant\Bbb O$. Such a
factorization is not unique: if an operator $\Phi$ is triangular
and satisfies $\Phi^*\Phi=\Bbb I$ then $\Phi V$ also provides a TF
to $C$. Below a {\it canonical} factor $V$ is singled out.

\subsubsection*{Diagonal}

Let $\mathscr H$ be one more Hilbert space. For an operator
$W:\mathscr F\to\mathscr H$, the subspaces $\overline{W\mathscr
F_s},\,\,s\in[0,T]$ form a nest\, (but not necessarily continuous
\cite{BSim_AA_2024}). By $P_s$ we denote the projections in
$\mathscr H$ onto $\overline{W\mathscr F_s}$. Fix a $T<\infty$ and
choose a partition $\Xi=\{s_k\}_{k=0}^n:\quad
0=s_0<s_1<\dots<s_n=T $ of $[0,T]$ of the range
$r^{\Xi}_W:=\max\limits_{k=1,\dots,K}(s_k-s_{k-1})$. Denote
$\Delta X_{s_k}:=X_{s_k}-X_{s_{k-1}}$, $\Delta
P_{s_k}:=P_{s_k}-P_{s_{k-1}}$, and put
\begin{equation*}
D^{\Xi}_W\,:=\,\sum\limits_{k=1}^n\Delta P_{s_k}\,W\,\Delta
X_{s_k}\,.
\end{equation*}
\begin{Definition}\label{D Diag}
The operator
\begin{equation*}
D_W\,:=\,{\rm w}\text{-}\!\lim\limits_{r^\Xi_W\to 0}\,D^{\Xi}_W
=\int_{[0,T]}dP_s\,W\,dX_s
 \end{equation*}
is called a diagonal of $W$ w.r.t. the nest $\mathfrak f$.
\end{Definition}

\noindent The limit is understood by Riemann: for arbitrary
$\varepsilon>0$ and elements $f\in{\mathscr F},\,h\in{\mathscr H}$
there is $\delta>0$ such that
\begin{equation*}
|\,([D_W- D^{\Xi}_W]f,h)\,|\,<\,\varepsilon
\end{equation*}
holds for any partition $\Xi$ of the range $r^{\Xi}<\delta$. Such
a limit does not always exist (A.B.Pushnitskii \cite{BPush},
\cite{Dav}).

The diagonal intertwines the nests: the relation
\begin{equation}\label{Eq PW=WX}
D_WX_s\,=\,P_sD_W,\quad D^*_WP_s\,=\,X_sD^*_W  \qquad s\in[0,T]
 \end{equation}
holds. For its adjoint one has $D_W^*=\int_{[0,T]}dX_sW^*dP_s$,
where the integral converges (or diverges) in the same sense as
for $D_W$. The relation $\|D_W\|\leqslant\|W\|$ is valid
\cite{BSim_AA_2024}.

\subsubsection*{Factorization}

In what follows we deal with $\mathscr H=\mathscr F$. Let $C$ be a
positive operator in $\mathscr F$. Assume that its positive square
root $\sqrt C$ has the diagonal $ D_{\sqrt
C}=\int_{[0,T]}dP_s\sqrt C\,dX_s$, where $P_s$ projects in
$\mathscr F$ onto $\overline{\sqrt{C}\mathscr F_s}$. As is easily
seen from (\ref{Eq PW=WX}), the operator $V:=D^*_{\sqrt C}\,{\sqrt
C}$ is triangular: $V\mathscr F_s\subset\mathscr F_s$ holds. If,
in addition, the diagonal satisfies the conditions
\begin{equation}\label{Eq Isom Diag}
{\rm Ran\,} D_{\sqrt C}=\mathscr F,\qquad    D_{\sqrt
C}\,D^*_{\sqrt C}=\Bbb I,
\end{equation}
then the TF
\begin{equation}\label{Eq CTF}
C=V^*V,\qquad V=D^*_{\sqrt C}\,{\sqrt C}
\end{equation}
holds and is referred to as {\it canonical}. Its peculiarity and
advantage is that the triangular factor $V$ is determined
constructively via the factorizable operator $C$
\cite{BSim_AA_2024}. The conditions (\ref{Eq Isom Diag}) are
motivated by applications: they are realized in multidimensional
inverse problems \cite{BKach Oper Integral 94,B Obzor IP 97,B
IPI}.

\subsubsection*{Stability}

Let us prove Theorem \ref{T0}.
\begin{proof}
Fix $\varepsilon>0$ and arbitrary $f,g\in\mathscr F$.
Representing:
\begin{align*}
& V-V^\alpha=D^*_{\sqrt C}{\sqrt
C}-D^*_{\sqrt{C^\alpha}}\sqrt{C^\alpha}=(D^*_{\sqrt
C}-D^{\Xi\,*}_{\sqrt C})\sqrt
C-(D^*_{\sqrt{C^\alpha}}-D^{\Xi\,*}_{\sqrt{C^\alpha}})\sqrt{C^\alpha}+\\
& +
(D^{\Xi\,*}_{\sqrt{C}}-D^{\Xi\,*}_{\sqrt{C^\alpha}})\sqrt{C}+D^{\Xi\,*}_{\sqrt{C^\alpha}}(\sqrt{C}-\sqrt{C^\alpha}),
\end{align*}
we have
\begin{align*}
& |([V-V^\alpha]f,g)|\leqslant  \big|\left(\sqrt{C}f,(D_{\sqrt
C}-D^{\Xi}_{\sqrt C})g\right)\big|+\big|\left(\sqrt{C^\alpha}f,(D_{\sqrt{C^\alpha}}-D^{\Xi}_{\sqrt{C^\alpha}})g\right)\big|+\\
&+\big|\left(\sqrt{C}f,(D^{\Xi}_{\sqrt{C}}-D^{\Xi}_{\sqrt{C^\alpha}})g\right)\big|+\big|\left((\sqrt{C}-\sqrt{C^\alpha})f,D^{\Xi}_{\sqrt{C^\alpha}}g\right)\big|\\
& =:I+I\!I+I\!I\!I+I\!V.
\end{align*}
Taking $r^\Xi<r$ for a small enough $r$, we provide
$I<\frac{\varepsilon}{4}$ and $I\!I<\frac{\varepsilon}{4}$
uniformly w.r.t. $\alpha$. In $I\!I\!I$ one has
\begin{align*}
&
(D^{\Xi}_{\sqrt{C}}-D^{\Xi}_{\sqrt{C^\alpha}})\,g=\sum\limits_{k=1}^N\left[\Delta
P_{s_k}\sqrt{C}-\Delta P^\alpha_{s_k}\sqrt{C^\alpha}\right]\Delta
X_{s_k}g=\\
&=\sum\limits_{k=1}^N\left[\left(P_{s_k}\sqrt{C}-P^\alpha_{s_k}\sqrt{C^\alpha}\right)-\left(P_{s_{k-1}}\sqrt{C}-P^\alpha_{s_{k-1}}\sqrt{C^\alpha}\right)\right]\widetilde
g_k,
\end{align*}
where $\widetilde g_k:=\Delta X_{s_k}g=X_{s_k}\Delta X_{s_k}g\in
X_{s_k}\mathscr F$. By the latter embedding,
$$
P_{s_k}\sqrt{C}\,\widetilde g_k=\sqrt{C}\,\widetilde
g_k\quad\text{and}\quad P^\alpha_{s_k}\sqrt{C^\alpha}\,\widetilde
g_k=\sqrt{C^\alpha}\,\widetilde g_k.
$$
holds and implies
\begin{align*}
&
I\!I\!I\leqslant\sum\limits_{k=1}^N\left[\|(\sqrt{C}-\sqrt{C^\alpha})\,\widetilde
g_k\|+
\|(P_{s_{k-1}}\sqrt{C}-P^\alpha_{s_{k-1}}\sqrt{C^\alpha})\,\widetilde
g_k\|\right]=\\
&=\sum\limits_{k=1}^N\left[\|(\sqrt{C}-\sqrt{C^\alpha})\,\widetilde
g_k\|+
\|(P_{s_{k-1}}-P_{s_{k-1}}^\alpha)\sqrt{C}+P^\alpha_{s_{k-1}}(\sqrt{C}-\sqrt{C^\alpha})\,\widetilde
g_k\|\right]\leqslant\\
& \leqslant
\sum\limits_{k=1}^N\left[2\,\|(\sqrt{C}-\sqrt{C^\alpha})\,\widetilde
g_k\|+ \|(P_{s_{k-1}}-P_{s_{k-1}}^\alpha)\sqrt{C}\,\widetilde
g_k\|\right].
\end{align*}
Since $\sqrt{C^\alpha}\overset{\rm reg}\to\sqrt{C}$ means
$\sqrt{C^\alpha}\overset{\rm s}\to\sqrt{C}$ and
$P^\alpha_{s_{k-1}}\overset{\rm s}\to P_{s_{k-1}}$, for big enough
$\alpha$ we provide $I\!I\!I,I\!V<\frac{\varepsilon}{4}$ and
eventually arrive at $|([V-V^\alpha]f,g)|\leqslant \varepsilon$.
Thus, the convergence $V^\alpha\overset{\rm w}\to V$ does occur.
\end{proof}

\subsubsection*{Comments and examples}

\noindent$\bullet$\,\,\, The notion of the regular convergence
$\overset{\rm reg}\to$ is motivated by the following question of
general character. Let $\mathscr H$ be a Hilbert space, $\mathscr
M\subset\mathscr H$ a subspace. Assume that $W^\alpha\to W$ holds
in a certain operator topology. Let $P^\alpha$ and $P$ be the
projections onto $\overline{W^\alpha\mathscr M}$ and
$\overline{W\mathscr M}$ respectively. Does the convergence
$P^\alpha\to P$ necessarily hold? Evidently, the answer is
negative: for instance, $\alpha^{-1}\Bbb I=:W^\alpha\to W=\Bbb O$
by norm but, taking $\mathscr M=\mathscr H$, we have $\mathbb
I=P^\alpha\not\to P=\Bbb O$. The example below shows that the
positivity of $W^\alpha$ and $W$ does not help. In such a
situation, the sufficient (verifiable in applications) conditions
on $W^\alpha$, $W$ and $\mathscr M$ which ensure $P^\alpha\to P$
are important and interesting.
\medskip

\noindent$\bullet$\,\,\, There is an example. Let
$\{\phi_k\}_{k=1}^\infty$, $(\phi_k,\phi_l)=\delta_{kl}$ be a
basis in $\mathscr H$, $\mathscr M:={\rm
span\,}\{\phi_k\}_{k=2}^\infty=\mathscr H\ominus{\rm
span\,}\{\phi_1\}$. Define $W$ and $W_n$ by
\begin{align*}
& W\phi_k\,=\,\frac{1}{k}\,\phi_k,\qquad k\geqslant 1
\end{align*}
and, for $n\geqslant 1$,
\begin{align*}
& W_n\phi_k=W\phi_k\,=\,\frac{1}{k}\,\phi_k,\qquad k\not=1,\,k\not= n,\\
& W_n\phi_1=\,\phi_1+\frac{1}{n}\,\phi_n,\\
& W_n\phi_n=\,\frac{1}{n}\,\phi_1+\frac{2}{n^2}\,\phi_n,\\
\end{align*}
so that the action of $W$ differs from the action of $W_n$ only on
the subspace ${\rm span\,}\{\phi_1,\phi_n\}$. Meanwhile, we have
$W^n\overset{\|\cdot\|}\to W$ as $n\to\infty$.

Take $\psi_n:=\phi_1-\frac{n}{2}\,\phi_n$,
$\|\psi_n\|^2=1+\frac{n^2}{4}$, and let
$\tilde\psi_n:=\|\psi_n\|^{-1}\psi_n$. Then, as is easy to check,
for the images of $\mathscr M$ and the projections on the images
we have
$$
W\mathscr M\,=\,\mathscr M,\qquad P=\mathbb I-(\cdot,\phi_1)\phi_1
$$
and
$$
W_n\mathscr M\,=\,\mathscr H\ominus {\rm span\,}\{\psi_n\},\qquad
P_n=\mathbb I-(\cdot,\tilde\psi_n)\tilde\psi_n\,.
$$
Tending $n\to\infty$, we get $P_n\overset{\rm s}\to\Bbb I\not=P$.
\medskip

\noindent$\bullet$\,\,\, Let $C$ and $C^\alpha$ be {\it positive
definite} operators (isomorphisms) admitting the canonical TF, and
let the convergence $C^\alpha\to C$ holds by the operator norm.
Then $\sqrt{C^\alpha}\overset{\rm reg}\to\sqrt{C}$ holds. Indeed,
we have
\begin{align*}
P_s=\sqrt{C}\dot X_s\left(\dot X_s^*C\dot X_s\right)^{-1}\dot
X_s^*\sqrt{C}, \qquad s\in[0,T],
\end{align*}
where $\dot X_s$ is $X_s$ understood as an operator from $\mathscr
F$ to $\mathscr F_s$, so that $\dot X_s^*$ embeds $\mathscr F_s$
into $\mathscr F$ and the relations $\dot X_s\dot X_s^*=\Bbb
I_{\mathscr F_s}$, $\dot X_s^*\dot X_s=X_s$ hold. To verify this
representation it suffices to check the characteristic properties
of the projection $P_s$, which are $P_s^2=P_s$, $P_s^*=P_s$ and
${\rm Ran\,} P_s=\mathscr F_s$. Representing $P^\alpha_s$ in the
same way, one easily concludes that the strong convergence
$P^\alpha_s\to P_s$ occurs. So, the regularity of the convergence
$\sqrt{C^\alpha}\to\sqrt{C}$ on the nest $\mathfrak f$ does hold.

As one consequence of the above, in the case $C=\Bbb I+K>\Bbb O$
with a compact $K$, our canonical TF (\ref{Eq CTF}) coincides with
the classical Gokhberg-Krein special factorization and turns out
to be stable.
\medskip

\noindent$\bullet$\,\,\, One more case is as follows. Let $C$ and
$C^\alpha$,\,$\alpha>0$ admit the canonical TF (\ref{Eq 0}), and
let the integral sums defining $D_{\sqrt C}$ and $D_{\sqrt
{C^\alpha}}$ converge uniformly w.r.t. $\alpha$. Assume that
\begin{align*}
\mathscr F=\oplus\sum\limits_{l\geqslant 1}\mathscr F^l\,;\quad
F^lX_s=X_sF^l,\quad F^l\, C^\alpha=C^\alpha F^l, \quad
F^l\,C=CF^l,\qquad s\in[0,T]
\end{align*}
holds, where $F^l$ projects in $\mathscr F$ onto $\mathscr F^l$.
This means that $C^\alpha$ and $C$ (and hence $\sqrt{C^\alpha}$
and $\sqrt{C}$) are reduced by the subspaces ({\it channels})
$\mathscr F^l$.

Assume that $C^\alpha_l:=C^\alpha\!\upharpoonright\mathscr F_l$
and $C_l:=C\!\upharpoonright\mathscr F_l$ are the isomorphisms in
$\mathscr F_l$ for all $l$ and satisfy the conditions of Theorem
\ref{T0}. Then we have the TF in each channel:
$C^\alpha_l=V^{\alpha\,*}_lV^\alpha_l$, $C_l=V^*_lV_l$ and
$V^\alpha_l\overset{\rm w}\to V_l$ holds. As a result, we easily
get the factotisation $C^\alpha=V^{\alpha\,*}V^\alpha$, $C=V^*V$
and $V^\alpha\overset{\rm w}\to V$. In the mean time, if the low
bounds of the spectra of $C^\alpha_l$ and $C_l$ tend to zero as
$l\to \infty$, then neither $C^\alpha$ nor $C$ are isomorphisms.
\medskip

\noindent$\bullet$\,\,\, The latter case occurs in applications
but is not too rich in content. Omitting details, it is workable
in determination of the potential in the Schr\"odinger operator
$-\Delta+q$ in a ball with a spherically symmetric $q$. An
important open problem is to find at least sufficient conditions
that ensure $\sqrt{C^\alpha}\overset{\rm reg}\to\sqrt{C}$ while
being effectively verifiable in applications. As for the
determination of the potential, we mean the case of an arbitrary
domain and $q$ from a sufficiently representative class.

\noindent{\bf Key words:}\,\,triangular factorization, operator
diagonal, amplitude integral, canonical factorization, stability
of canonical factorization.
\smallskip

\noindent{\bf MSC:}\,\,\,47Axx,\,\, 47B25,\,\, 35R30.
\smallskip

\noindent{\bf UDK:}\,\,\,517.


\end{document}